\documentclass{amsart}

\usepackage{amsmath}
\usepackage{amsfonts}

\newcommand\nn{\mathbb N}
\newcommand\zz{\mathbb Z}
\newcommand\rr{\mathbb R}

\newcommand\cc{\mathbb C}
\newcommand\bb{\mathbb B}

\newcommand\oo{\mathcal O}
\def\epsilon{\varepsilon}
\def\phi{\varphi}
\newcommand{\twu}{\widetilde {\mathbb W}}
\newcommand{\wu}{\mathbb W}
\newcommand{\metrics}{\mathcal M}

\newcommand{\vol}{\operatorname{vol}}
\newcommand{\ds}{\Delta}
\newcommand{\conv}{\operatorname{conv}}

\theoremstyle{plain}
\newtheorem{theorem}{Theorem}
\newtheorem{proposition}[theorem]{Proposition}
\newtheorem{lemma}[theorem]{Lemma}
\newtheorem{corollary}[theorem]{Corollary}

\theoremstyle{remark}
\newtheorem{remark}[theorem]{Remark}

\numberwithin{equation}{section}
\numberwithin{theorem}{section}

\begin{document}

\title{On the Wu metric in unbounded domains}
\author{Piotr Jucha}
\address{Institute of Mathematics, Jagiellonian University, Reymonta 4,\newline 
30--059 Krak\'ow, Poland}
\email{Piotr.Jucha@im.uj.edu.pl}
\thanks{The work is supported by the Research Grant No.~1~PO3A~005~28 of the Polish Ministry 
of Science and Higher Education.}

\begin{abstract}
We discuss the properties of the Wu pseudometric
and present counterexamples for its upper semicontinuity that
answers the question posed by Jarnicki and Pflug.
We also give formulae for the Wu pseudometric in elementary Reinhardt domains.
\end{abstract}

\maketitle

\section{Introduction}

H.\ Wu introduced in \cite{bib:w} a new invariant metric which was to combine
invariant properties of the Kobayashi--Royden metric and regularity properties of 
K\"ahler metrics.
The metric depends on some initial (pseudo)metric $\eta$.
(Originally, it was defined only for the Kobayashi--Royden pseudometric.)
We call it the Wu (pseudo)metric associated to $\eta$ and write $\wu \eta$.

The pseudometric was studied in several
papers (e.g.\ \cite{bib:ck1}, \cite{bib:ck2}, \cite{bib:j1}, \cite{bib:j2},
 \cite{bib:jp3}, \cite{bib:jp2}).
Jarnicki and Pflug pointed out (\cite{bib:jp3}, \cite{bib:jp2}) that such an elementary property
as its upper semicontinuity had not been completely understood.
The question of the upper semicontinuity appears naturally,
for instance in the definition of the integrated form $\int(\wu\eta)$.
In general, the upper semicontinuity of $\eta$ does not imply the upper 
semicontinuity of $\wu \eta$ (cf.\ Remark~\ref{rem:one}).
If $\eta$ is the Kobayashi--Royden (pseudo)metric, the problem
has remained open, even though Wu (cf. \cite{bib:w_notes}, \cite{bib:w}) and
Cheung and Kim (cf. \cite{bib:ck1}) claimed (without proof) 
the upper semicontinuity of $\wu\kappa_D$. 

Jarnicki and Pflug asked then (\cite{bib:jp3}, \cite{bib:jp2})
whether $\wu \eta$ is upper semicontinuous if $\eta$ is one of the well-known pseudometrics:
Kobayashi--Royden ($\kappa$), Azukawa ($A$) or 
Carath\'eodory--Reiffen pseudometric of $k$--th order ($\gamma^{(k)}$).

We gave in \cite{bib:j2} an example of bounded pseudoconvex domain $D$ such that
the Wu metrics associated to $\kappa_D$ and $A_D$ are not
upper semicontinuous (cf.\ Proposition~\ref{ex:one}).
Moreover, it is known that if $D$ is a bounded domain then Carath\'eodory--Reiffen metrics
of any order are continuous (cf.\ \cite{bib:nik}). 
In view of Proposition~\ref{thm:wprop}~(\ref{en:cont}), the Wu metrics
associated to them are continuous as well.

We solve here the remaining unbounded case and give the full and negative answer
to Jarnicki and Pflug's question.
We would also like
to attract the attention to another aspect of the problem, which appears in unbounded domains.
Namely, the pseudometric $\wu$ is a normalization of the original metric introduced
in \cite{bib:w}---we denote it by $\twu$. For any admissible metric $\eta$
we have $\wu\eta(z;\cdot) = \sqrt{m(z)}\, \twu\eta(z;\cdot)$ where the constant $m(z)$
is the codimension of the subspace $\{X\in\cc^n:\, \widehat\eta(z;X)=0 \}$.
To justify the normalization, let us  mention a neat product formula 
(cf. Proposition~\ref{thm:wprop}~(\ref{en:prod})) but also an example of a domain
in which $\twu\kappa$ is not upper semicontinuous (cf.\ Remark~\ref{rem:two}).
This is why we also investigate the semicontiuity of $\twu$.
In a bounded domain the factor $m(z)$ does not depend on the point $z$, so it is 
irrelevant to the problem of semicontinuity.

We construct unbounded pseudoconvex Reinhardt domains $G_n\subset \cc^n$ for $n\ge 2$
such that for any contractible family
of pseudometrics $(\alpha_D)_{D\subset \cc^n}$ the pseudometrics
$\wu\alpha_{G_n}$ (for $n\ge 2$) and $\twu\alpha_{G_n}$ (for $n\ge 3$) are not
upper semicontinuous (the main results: Proposition~\ref{ex:two} and Proposition~\ref{ex:three}).
Moreover, we show (cf.\ Proposition~\ref{ex:monotone}) that $(\twu\alpha_D)_D$ and $(\wu\alpha_D)_D$ 
need not be monotone
(monotone here is understood as:
if $D_m\nearrow D$ then $\eta_{D_m} \to \eta_D$).
Recall that $(\alpha_D)_D$ is monotone
for $\alpha=\gamma^{(k)},A,\kappa$ (cf.\ \cite{bib:jp1}).

The above--mentioned results obviously give a negative answer to Jarnicki and Pflug's question.
Nevertheless, there is also a positive result which indicates that $2$--dimensional case is different.
Namely, 
the pseudometric $\twu\alpha_D$ is upper semicontinuous 
if $\alpha_D$, for a $2$--dimensional domain $D$, is a continuous pseudometric
(cf. Proposition~\ref{prop:cont}).
In particular, $\twu\gamma_D$ is upper semicontinuous if $D\subset \cc^2$.
We do not know whether the same is true for
Carath\'eodory--Reiffen metrics of higher order.

Since most considerations involve only unbounded domains, there appears a question 
(suggested by Professor M.~Jarnicki): 
Is there an $\eta$--hyperbolic (or pointwise $\eta$--hyperbolic)
domain $D$ such that $\wu\eta_D$ is not upper semicontinuous?
Certainly, Proposition~\ref{ex:one} gives the answer for Kobayashi--Royden and Azukawa metrics.
The problem remains open for $\eta=\gamma^{(k)}$. However, such a domain does not
exist in the class of pseudoconvex Reinhardt domains (cf.\ Proposition~\ref{prop:reinh}).

In the last section we present the formulae for the Wu metric in elementary Reinhardt domains
(cf. Proposition~\ref{prop:formulae}).
The formula for $\twu\kappa$ has been already given in \cite{bib:j1}.

\section{Definition and known facts}

We denote by $\ds$ the open unit disk in $\cc$.
Let $(\eta_D)_D$ be a family of pseudometrics
defined for all domains $D\subset \cc^n$, $n\ge 1$, i.e.
\[
\eta_D: \, D\times \cc^n \to \rr_+, \quad \eta(a; \lambda X) = |\lambda| \eta(a;X), \quad
\lambda \in \cc, (a,X) \in D \times \cc^n.
\label{eq:mdef1}
\]

We call $(\eta_D)_D$ a \emph{holomorphically contractible family of pseudometrics}
if the following two conditions are satisfied
(cf. \cite{bib:jp1}, \cite{bib:jp2}):
\begin{align*}
\eta_\ds(z;X) &= \frac {|X|}{1-|z|^2}, \qquad z\in \ds, X \in \cc,\\
\eta_{D_2}(F(z);F^\prime(z)X) &\le \eta_{D_1}(z;X), \quad z\in D_1, X\in \cc^{n_1},
\end{align*}
for any domains $D_1\subset \cc^{n_1}$, $D_2 \subset \cc^{n_2}$ and every 
holomorphic mapping $F:D_1\to D_2$.

We say that the family $(\eta_D)_D$ has \emph{the product property} if
\begin{multline}
\eta_{D_1\times D_2}((z,w);(X,Y)) = \max (\eta_{D_1}(z;X), \eta_{D_2}(w;Y)),\\
(z,w)\in D_1\times D_2, (X,Y) \in\cc^{n_1}\times \cc^{n_2}
\label{eq:prod}
\end{multline}
for any domains $D_1\subset \cc^{n_1}$, $D_2\subset \cc^{n_2}$, $n_1,n_2 \ge 1$.

A domain $D\subset \cc^n$ is called \emph{$\eta$--hyperbolic} if
\begin{equation}
\forall\, a \in D\ \ \exists\, C,r>0\ \ \forall\, z\in\bb(a,r)\cap D, X\in\cc^n: 
  \eta_D(z;X)\ge C\|X\|.
\label{eq:hyper}
\end{equation}
A domain $D\subset \cc^n$ is \emph{pointwise $\eta$--hyperbolic} if
\[
\forall\, z\in D, X\in\cc^n\setminus\{0\}: \eta_D(z;X)>0.
\]
In the latter case we call $\eta_D$ a \emph{metric}.

Recall definitions of the $k$--th order Carath\'eodory--Reiffen ($\gamma^{(k)}$),
Azukawa ($A$) and Kobayashi--Royden ($\kappa$) pseudometrics.
For details and properties see e.g.\ \cite{bib:jp1}, \cite{bib:jp2}.

For a domain $D \subset \cc^n$ and $a \in D$, $X\in\cc^n$, $k\in\nn\setminus\{0\}$ define:
\begin{align*}
\gamma_D^{(k)}(a;X)&:= \sup\Big\{\Big|\frac{1}{k!}f^{(k)}(a)X \Big|^{\frac 1k}:\, 
  f\in \oo(D,\ds), \textrm{ ord}_a f\ge k \Big\};\\
A_D(a;X)&:=\sup\Big\{ \limsup_{0\neq \lambda \to 0} \frac {v(a+\lambda X)}{|\lambda|}:\,
  v:D\to [0,1) \textrm{ is log--psh},\\
    &\qquad\qquad \exists\, M,r>0:\, v(z)\le M\|z-a\| \textrm{ if } \|z-a\| < r \Big\};\\
\kappa_D(a;X)&:= \inf\{t>0:\, \exists\, \phi\in\oo(\ds,D):\, 
  \phi(0)=a,\, t\phi^\prime(0) =X \}.
\end{align*}
We write $\gamma_D:= \gamma_D^{(1)}$.

\smallskip

We present the sketch of the definition of the Wu metric in an abstract setting (\cite{bib:jp3}).
For detailed discussion we refer the reader to \cite{bib:jp3} (or \cite{bib:jp2}) 
and \cite{bib:w}.

For a domain $D \subset \cc^n$, denote by $\metrics (D)$ the space of all pseudometrics
such that
\begin{equation}
\forall\, a \in D\ \exists\, M,r>0 :\, \eta(z;X) \le M\|X\|,
  \ z \in \bb(a,r) \subset D, X \in \cc^n,
\label{eq:mdef2}
\end{equation}
where $\bb(a,r):=\{z\in \cc^n:\, \|z-a\| <r \}$ and $\|\cdot\|$ is the standard Euclidean norm. 
Note that condition (\ref{eq:mdef2}) is satisfied if $\eta$ is upper semicontinuous.

For convenience, let
$\bb_{\eta_D}(a):= \{X\in\cc^n:\, \eta_D(a;X)<1\}$ 
be the unit $\eta_D$--ball at a point $a\in D$
(for $\eta_D\in\metrics (D)$).

Let $\widehat\eta$ denote the Busemann pseudometric 
associated to $\eta$ (cf.\ e.g.\ \cite{bib:jp1}), i.e.
\[
\widehat\eta(a;X) := \sup\{p(X)\}, \quad a \in D, X\in\cc^n,
\]
where the supremum is taken over all $\cc$--seminorms $p$ such that
$p\le \eta(a;\cdot)$.
We have $\widehat\eta\le \eta$. Recall that if $\eta$ is upper semicontinuous,
then so is $\widehat\eta$ and
 $\bb_{\widehat\eta}(a)=\conv \bb_\eta(a)$.

Fix a domain $D\subset \cc^n$, a point $a\in D$, a pseudometric $\eta \in \mathcal M(D)$ and put:
\begin{align*}
V_{\eta}(a):= &\{X\in\cc^n:\, \widehat\eta(a;X)=0\},\\
U_{\eta}(a):= &\textrm{ the orthogonal complement of $V_{\eta}(a)$ with respect to}\\
& \textrm{ \hfill the standard scalar product
in } \cc^n.
\end{align*}

For any pseudo--Hermitian scalar product $s:\cc^n\times\cc^n \to \cc$, define
\begin{align*}
q_s(X):=&\sqrt{s(X,X)}, \quad X\in \cc^n.
\end{align*}

Let $\mathcal F(\eta,a)$ be a set of all pseudo--Hermitian scalar products 
$s:\cc^n\times\cc^n \to \cc$
such that $q_s\le \eta(a;\cdot)$ (or, equivalently, $\bb_{\eta}(a)\subset \bb_{q_s}$).
There exists a unique (!) element $s(\eta,a) \in \mathcal F(\eta,a)$ 
that is maximal with respect to
the partial ordering $\prec$ defined for $\alpha,\beta\in\mathcal F(\eta,a)$:
\[
\alpha \prec \beta\ \textrm{ if } \ \det[\alpha(e_j,e_k)]_{j,k=1,\dots,m}
    \le \det[\beta(e_j,e_k)]_{j,k=1,\dots,m},
\]
for any basis $(e_1,\dots,e_m)$ of $U_{\eta}(a)$.

We define
\begin{align}
\twu \eta(a;X) &:= q_{s(\eta,a)}(X), \quad X \in \cc^n;\notag\\
\wu \eta(a;X) &:= \sqrt {m(\eta,a)} \, \twu \eta(a;X), \quad X \in \cc^n,
\label{eq:wu}
\end{align}
where $m(\eta,a):= \textrm{dim}\, U_{\eta}(a)$. 

Note that the definition of $\wu\eta$ depends, in fact, only on $\widehat\eta$.
Moreover, the construction determines that
the ball $\bb_{\twu\eta}(a)$ is the ``minimal'' ellipsoid 
containing $\bb_{\eta}(a)$ (in volume if $\bb_{\eta}(a)$ is bounded).

Some basic properties of the Wu pseudometric are listed in Proposition~{\ref{thm:wprop}}.

\begin{proposition}[cf.\ \cite{bib:w}, \cite{bib:jp3}, \cite{bib:jp2}]{\ }\label{thm:wprop}
\begin{enumerate}
\item\label{en:cont}
If $\eta\in\metrics(D)$ is a continuous complete metric, then so is $\wu\eta$.

\item If $(\eta_D)_{D}$ 
  is a holomorphically contractible family of pseudometrics,
then for any  biholomorphic mapping $F:D_1\to D_2$  $(D_j \subset \cc^n$, $j=1,2)$ we have
\[
\wu \eta_{D_2}(F(z);F^\prime(z)X) = \wu\eta_{D_1}(z;X), \quad z\in D_1, X\in \cc^n.
\]

\item
If $(\eta_D)_{D}$ 
is a holomorphically contractible family of pseudometrics,
then for any  holomorphic mapping $F:D_1\to D_2$  $(D_1 \subset \cc^{n_1}$, 
$D_2 \subset \cc^{n_2})$ we have
\[
\wu\eta_{D_2}(F(z);F^\prime(z)X) \le \sqrt{n_2} \wu\eta_{D_1}(z;X), \quad z\in D_1, X\in \cc^{n_1}.
\]

\item \label{en:prod}
If $(\eta_D)_D$ is a family of pseudometrics satisfying  the product property, then
\begin{multline*}
\wu\eta_{D_1\times D_2}((z,w);(X,Y)) 
  = \big( \big(\wu \eta_{D_1}(z;X)\big)^2 + \big(\wu \eta_{D_2}(w;Y)\big)^2\big)^{\frac 12},\\
(z,w)\in D_1\times D_2 \subset \cc^{n_1}\times \cc^{n_1}, (X,Y) \in\cc^{n_1}\times \cc^{n_2},
\end{multline*}
\end{enumerate}
\end{proposition}

\smallskip

In the next section we shall use Lemma~\ref{prop:c2}, which comprises
some of the properties of $\twu\eta$--balls. Its proof in the two--dimensional case
is essentially contained in \cite{bib:j2}.

Consider the following mapping (cf. \cite{bib:ck1}, \cite{bib:ck2})
\[
\Psi:\, \cc^n \to \rr_+^n, \quad \Psi(z):= (|z_1|^2,\dots,|z_n|^2), \quad z\in\cc^n.
\]
Note that $\Psi$ transforms any bounded complete Reinhardt ellipsoid in $\cc^n$ into
a simplex
\[
T_a:= \big\{ (u_1,\dots,u_n)\in \rr_+^n:\, \sum_{j=1}^n \frac {u_j}{a_j} <1\big\}
\]
for some $a=(a_1,\dots,a_n)\in (\rr_+\setminus\{0\})^n$.
In fact,
the mapping $\Psi$ determines the one-to-one correspondence between bounded complete Reinhardt 
ellipsoids in $\cc^n$ and simplexes $T_{a}\subset \rr_+^n$.


\begin{lemma}[cf.\ \cite{bib:j2}] \label{prop:c2}
Let $D$ be a domain in $\cc^n$, $z_0\in D$, and $\eta\in\metrics(D)$ be a pseudometric such that
$\bb_{\eta}(z_0)$ is a bounded Reinhardt domain. Then
\begin{enumerate}
\item\label{en:c2_1} $\bb_{\twu\eta}(z_0)$ is a complete Reinhardt domain;
\item\label{en:c2_2} there exist $a_1,\dots,a_n>0$ such that 
  $\twu\eta(z_0;X)= \big( \sum_{j=1}^n\frac {|X_j|^2}{a_j} \big)^{\frac 12}$
  for $X=(X_1,\dots,X_n)\in\cc^n$;
\item\label{en:c2_3} 
$\Psi(\bb_{\twu\eta}(z_0))=T_{(a_1,\dots,a_n)}$ is a unique simplex of smallest volume 
that contains
  $\Psi(\bb_\eta(z_0))$, where the numbers $a_j$ are as in {\rm (\ref{en:c2_2})};
\item\label{en:c2_4}
if $\bb_{\eta}(z_0)=r_1 \Delta\times\dots\times r_n \Delta$ for some $r_1,\dots,r_n>0$, 
then $\Psi(\bb_{\twu\eta}(z_0)) = T_{nr_1^2,\dots,nr_n^2}$.
\end{enumerate}
\end{lemma}

\begin{proof}
(\ref{en:c2_1})
Since the ball $\bb_{\eta}(z_0)$ is invariant
under the action of the (volume preserving) transformations
\[
\Phi_{\Lambda}(X):=(\lambda_1X_1,\dots,\lambda_nX_n), \quad X\in\cc^n,
\Lambda=(\lambda_1,\dots,\lambda_n)\in (\partial \ds)^n,
\]
then so is 
$\bb_{\twu\eta}(z_0)$. Otherwise, it would contradict its uniqueness.
Moreover,  $\bb_{\twu\eta}(z_0)$ is convex, and consequently complete Reinhardt.

(\ref{en:c2_2})
Let $[a_{jk}]_{j,k=1,\dots,n}$ be the matrix representation of
the Hermitian scalar product associated with $\twu\eta(z_0;\cdot)$
in the canonical basis of $\cc^n$, i.e.\
$\twu\eta(z_0;X)^2 = \sum_{j,k=1}^n a_{jk}X_j\bar X_k$, $X\in\cc^n$.
The invariance of $\bb_{\twu\eta}(z_0)$ under $\Phi_\Lambda$ implies that
$a_{jk}=0$ for $j\neq k$.
Certainly, $a_j=a_{jj}^{-1}$.

(\ref{en:c2_3})
Note that volumes of complete Reinhardt ellipsoids and 
the corresponding triangles $T_a$ are proportional. Namely,
$\vol \{(X_1,\dots,X_n)\in\cc^n :\, \sum_{j=1}^n \frac{X_j}{a_j} <1\} = n\beta_n \vol T_{(a_1,\dots,a_n)}$,
where $\beta_n$ denotes the volume of the Euclidean unit ball in $\cc^n$.

(\ref{en:c2_4})
It suffices to minimize the volume $V(b):= \vol T_b$ in the set
$\{b\in\rr_+^n:\, \sum_{j=1}^n \frac {r_j^2}{b_j} \le 1 \}$.
\end{proof}

\begin{remark}[cf.\ \cite{bib:jp2}]\label{rem:one}
The main reason why the Wu metric is not always upper semicontinuous is 
illustrated by the following example.
Let $D\subset \cc^2$ be any domain, $z_0\in D$,  and 
$\alpha : D \times \cc^2\to \rr_+$ a pseudometric defined as follows:
\[
\alpha(z;(X_1,X_2)) := 
  \begin{cases}
    |X_1|^2 +|X_2|^2, &\text{ if } z\neq z_0\\
    \max\{|X_1|, \frac 12 |X_2|\}, &\text{ if } z=z_0
  \end{cases},
  \quad \text{ } (X_1,X_2)\in\cc^2.
\]
Although $\bb_{\alpha}(z) \subset \bb_{\alpha}(z_0)$ for all $z\in D\setminus\{z_0\}$,
there is no inclusion between $\bb_{\wu\alpha}(z)$ and $\bb_{\wu\alpha}(z_0)$.
Hence, $\alpha$ is upper semicontinuous but $\wu\alpha$ is not.
\end{remark}

\begin{remark}[cf.\ \cite{bib:jp2}]\label{rem:two}
Likewise, the pseudometric $\twu\kappa$ is not always upper semicontinuous.
Let $D\subset \cc^2$ and $D\ni z_k\to z_0 \in D$ be such that $\kappa_D(z_k;\cdot)$ is not a metric
and $\kappa_D(z_0;\cdot)$ is a metric. (Such a domain exists, cf.\ \cite{bib:jp1}.)
Put $G:=D\times \ds\subset \cc^3$. Then
\[
\twu\kappa_D((z_k,0);(0,0,1)) \ge \frac 1{\sqrt 2},\quad
\twu\kappa_D((z_0,0);(0,0,1)) = \frac 1{\sqrt 3}.
\]
Therefore, $\twu\kappa_D$ is not upper semicontinuous.
\end{remark}

\section{Results on upper semicontinuity}

We assume in the sequel that the family $(\alpha_D)_D$ defined for all domains $D\subset\cc^n$, $n\ge 1$,
is a holomorphically contractible family of pseudometrics.
Consequently, we have that $\gamma_D\le \alpha_D\le \kappa_D$ 
(cf.\ \cite{bib:jp1}) and $\alpha_D\in\metrics(D)$ for any domain $D\subset \cc^n$.

\begin{proposition}[\cite{bib:j2}]\label{ex:one}
Define
\[
G:= \{(z_1,z_2)\in\cc^2:\, |z_1|<1, |z_2|<5, 10|z_2|e^{u(z_1)}<1\},
\]
where
$u(z_1) = 1 + \sum_{j=4}^\infty \frac 1{2^j} \max \{ \log \frac {|2^{-j}-z_1|}{2}, -2^{2j}\}$.
If $A_{D}\le \alpha_D\le \kappa_{D}$ for all domains $D\subset \cc^2$, 
then neither $\wu\alpha_G$ nor $\twu\alpha_G$ is upper semicontinuous.
\end{proposition}

\begin{remark}
A similar domain as in Proposition~\ref{ex:one}
can be constructed in higher dimensions
if $(\alpha_D)_D$ (for $D\subset\cc^n, n\ge 1$) is a family of pseudometrics satisfying the product
property (\ref{eq:prod}). 
It suffices to take the Cartesian product $G\times \ds^{n-2}$
and use Proposition~\ref{thm:wprop}~(\ref{en:prod}).
\end{remark}

\begin{proposition}\label{ex:two}
Define 
\[
G_2:= \{ (z_1,z_2) \in \cc^2:\, |z_1|(1+|z_2|) < 1\}.
\]
Then $\wu\alpha_{G_2}$ is not upper semicontinuous.
\end{proposition}
\begin{proof}
We shall prove that
\begin{equation}
\begin{split}
\limsup_{0<x\to0} \wu\alpha_{G_2}\big((x,0);(1,0)\big)
& = \sqrt 2 \limsup_{0<x\to0} \twu\alpha_{G_2}\big((x,0);(1,0)\big) \\
& \ge \sqrt 2  > 1 = \wu\alpha_{G_2}\big((0,0);(1,0)\big).
\end{split}
\label{eq:cond2}
\end{equation}
Above, the factor $m(\alpha_{G_2})$ is crucial (cf.~(\ref{eq:wu})).
We have $m(\alpha_{G_2},(x,0))=2$ for $x\in(0,1)$ because of boundedness of 
$\bb_{\twu\alpha_{G_2}}(x,0)$ (cf.\ \emph{Step~1}) while $m(\alpha_{G_2},(0,0))=1$.

\emph{Step~1. The balls $\bb_{\alpha_{G_2}}(x,0)$ for $x\in(0,1)$ are bounded
Reinhardt domains.}

First, note that they
are Reinhardt domains.
Indeed, rotations of the form
$\cc^2 \ni (X_1,X_2)\mapsto (X_1,\lambda X_2)$
are automorphisms of $G_2$ and they fix points $(x,0)$.
Due to the contractibility of $\alpha$, the balls $\bb_{\alpha_{G_2}}(x,0)$ are
also invariant under these rotations and, moreover, they are balanced.

Now, take the mapping $F(z_1,z_2):= z_1(1+z_2)$ for $(z_1,z_2)\in \cc^2$. 
We have $F(G_2) \subset \ds$ and for $X_1,X_2>0$ we get
\[
\gamma_{G_2} \big((x,0); (X_1,X_2) \big) \ge \gamma_\ds \big( F(x,0); F^\prime(x,0)(X_1,X_2) \big)
= \frac {X_1 +xX_2}{1-x^2}.
\]
Since the balls $\bb_{\alpha_{G_2}}(x,0)$ are Reinhardt,
we have also
\begin{equation}
\bb_{\alpha_{G_2}}(x,0) \subset \{(X_1,X_2)\in\cc^2:\, \frac {|X_1| +x|X_2|}{1-x^2}<1 \},
\quad x \in (0,1)
\label{eq:ball}
\end{equation}
which implies the boundedness of $\bb_{\alpha_{G_2}}(x,0)$.


\emph{Step~2. $\twu\alpha_{G_2}((0,0);(X_1,X_2)) = \wu\alpha_{G_2}((0,0);(X_1,X_2))=|X_1|$,
for $(X_1,X_2)\in\cc^2$.}

From inequalities $\gamma_{G_2}\le \alpha_{G_2} \le \kappa_{G_2}$  we obtain
\[
\bb_{\gamma_{G_2}} (z) \supset \bb_{\alpha_{G_2}} (z) \supset \bb_{\kappa_{G_2}} (z),
  \quad z \in G_2.
\]
Since $G_2$ is a pseudconvex complete Reinhardt domain we have that $\bb_{\kappa_{G_2}}(0,0) = G_2$
and $\bb_{\gamma_{G_2}}(0,0)= \conv G_2 = \ds\times \cc$ (cf.\ \cite{bib:jp1}).
Hence, we get 
$\bb_{\twu\alpha_{G_2}}(0,0) 
=\ds\times\cc$
which implies the required formula.

\emph{Step~3. $\limsup_{0<x\to0} \twu\alpha_{G_2}((x,0);(1,0)) \ge 1$.}

Assume for a contradiction that there exist numbers $t>1$, $\delta\in(0,1)$ 
such that $\twu\alpha_{G_2}((x,0);(1,0)) < \frac 1t$ for any $x\in (0,\delta)$

Fix $x\in(0,\delta)$.
Since the ball $\bb_{\alpha_{G_2}}(x,0)$ is a bounded Reinhardt domain,
there exist numbers $a,b >0$ such that 
$T_{a,b}= \Psi(\bb_{\twu\alpha_{G_2}}(x,0))$
is the unique triangle of minimal area containing the set 
$\Psi(\bb_{\alpha_{G_2}}(x,0))$ (cf.\ Lemma~\ref{prop:c2}).
It follows from 
the assumption
that $(t,0)\in \bb_{\twu\alpha_{G_2}}(x,0)$, thus
$a>t^2$. On the other hand, we have $b> (\frac {1-x}{x})^2$ because
\begin{equation}
\alpha_{G_2} \Big( (x,0); \Big(0,\frac {1-x}{x}\Big)\Big) 
  \le \kappa_{G_2} \Big( (x,0); \Big(0,\frac {1-x}{x}\Big)\Big)
\le 1.
\label{eq:alpha}
\end{equation}
To get the latter inequality, take the function 
$\phi(\lambda):= (x, \frac {1-x}{x}\lambda)$, $\lambda\in\ds$.

Now, consider the triangle $T:= T_{(1,x^{-2})}$. 
Using condition~(\ref{eq:ball}) one can see that 
$\Psi(\bb_{\alpha_{G_2}}(x,0)) \subset T$.
We compare area of the both triangles:
\[
\frac {\vol T_{a,b}} {\vol T} = { x^2 a b}
> t^2(1-x)^2
\]
Therefore, for a sufficiently small $x$ we have $\vol T_{a,b} > \vol T$, which
contradicts the minimality of the triangle $T_{a,b}$.
\end{proof}

\begin{remark}\label{rem:g2}
It follows from the proof that for any $x\in(0,1)$
\begin{gather*}
\Big(0,\frac 1x -1\Big), (1-x^2,0) \in \overline{\bb_{\kappa_{G_2}}(x,0)}
\end{gather*}
Indeed, for the first point it is the direct consequence of (\ref{eq:alpha})---note that it is true
for all $x\in(0,1)$.
To verify the same for the point $(1-x^2,0)$, take the mapping
$\phi(\lambda)= (\frac{\lambda+x}{1+x\lambda},0)$, $\lambda \in\ds$.
\end{remark}

\begin{proposition}\label{ex:three}
For $n\ge 3$ define
\[
G_n:= G_2\times \ds^{n-2} \subset \cc^n.
\]
Then neither $\twu\alpha_{G_n}$ nor $\wu\alpha_{G_n}$ 
is upper semicontinuous.
\end{proposition}

\begin{proof}
We shall proceed in much the same way as in the proof of Proposition~\ref{ex:two}.
We are going to show that
\[
\begin{split}
\limsup_{z\to 0} \twu\alpha_{G_n}(z;(1,0,\dots,0))& \ge \sqrt\frac 2n
> \frac1{\sqrt{n-1}} = \twu\alpha_{G_n}(0;(1,0,\dots,0)),\\
\limsup_{z\to 0} \wu\alpha_{G_n}(z;(1,0,\dots,0))&\ge \sqrt{2}
>1 = \wu\alpha_{G_n}(0;(1,0,\dots,0).
\end{split}
\]

\emph{Step~1. The balls $\bb_{\alpha_{G_n}}(x,0,\dots,0)$ for $x\in(0,1)$ are bounded
Reinhardt domains.}

Recall that both families $(\gamma_D)_D$ and $(\kappa_D)_D$ satisfy the product 
property~(\ref{eq:prod}) and $\gamma_D\le \alpha_D\le \kappa_D$.
Hence, we have
\begin{equation}
\bb_{\kappa_{G_2}}(z) \times \ds^{n-2} \subset \bb_{\alpha_{G_n}}(z,w) 
  \subset \bb_{\gamma_{G_2}}(z) \times \ds^{n-2}, 
\ (z,w)\in G_2\times\ds^{n-2}.
\label{eq:balls}
\end{equation}
Therefore, the balls $\bb_{\alpha_{G_n}}(x,0,\dots,0)$ are bounded for $x> 0$.
They are also Reinhardt domains---note that they are balanced and invariant under 
rotations $\cc^n\ni X\mapsto (X_1,\lambda_2X_2,\dots,\lambda_nX_n)$ for 
$\lambda_j\in\partial\ds$, $j=2,\dots,n$.

\emph{Step~2. The following formula holds:}
\begin{equation}
\sqrt{n-1}\,\twu\alpha_{G_n}(0;(1,0,\dots,0))
  = \wu\alpha_{G_n}(0;(1,0,\dots,0)) = 1.
\label{eq:wugn}
\end{equation}

Recall that $\bb_{\kappa_{G_2}}(0)=G_2$.
Hence, we have that 
$\conv \bb_{\alpha_{G_n}}(0) = \bb_{\gamma_{G_n}}(0)= \ds\times\cc\times\ds^{n-2}$
by condition (\ref{eq:balls}), and consequently we get the required formula.

\emph{Step~3.
$\limsup_{0<x\to0}\twu\alpha_{G_n}((x,0,\dots,0);(1,0,\dots,0) \ge \sqrt{\frac 2n}$.}

Assume the contrary, i.e.\ there exist numbers $t>{\frac n2}$ and $\delta>0$ such that
$\twu\alpha_{G_n}((x,0,\dots,0);(\sqrt t,0,\dots,0)) <1$ for any $x\in(0,\delta)$.

Fix such an $x$. 
Since $\bb_{\alpha_{G_n}}(x,0,\dots,0)$ is a bounded Reinhardt domain
there exist an $n$--tuple $a=(a_{1},\dots,a_{n})\in(\rr_+\setminus\{0\})^n$ such that
$\Psi(\bb_{\twu\alpha_{G_n}}(x,0,\dots,0))=T_{a}$
(cf.~Lemma~\ref{prop:c2}).
Recall that $T_{a}$ has smallest volume of all simplexes containing the set
$\Psi(\bb_{\alpha_{G_n}}(x,0,\dots,0))$.
The assumption is then equivalent to inequality $a_{1}>t$.

Let $T:=T_{(\frac n2, \frac n{2x^2}, n,\dots,n)}$ be another simplex. Note that
\[
T \supset \Psi(\bb_{\gamma_{G_2}}(x,0) \times\ds^{n-2})
\supset \Psi\big(\bb_{\alpha_{G_n}}(x,0,\dots,0))
\]
because of condition~(\ref{eq:ball}) and (\ref{eq:balls}). 
From the minimality of $T_{a}$ we have that
$\vol T_{a}\le \vol T$. We shall estimate $\vol T_{a}$ and
show that the assumption $a_{1}>t>\frac n2$, in fact, leads to a 
contradiction, i.e.\ $\vol T_{a}> \vol T$ for small numbers $x>0$.

To simplify notation put $\mu:=(1-x^2)^2$, $\nu:=(\frac 1{x}-1)^2$.
From condition~(\ref{eq:balls}) and Remark~\ref{rem:g2} we obtain that
\[
(\mu,0,1,\dots,1),(0,\nu,1,\dots,1)\in \Psi(\overline{\bb_{\alpha_{G_n}}(x,0,\dots,0)})
 = \overline{T_{a}}.
\]
We shall find the simplex, say $T_{c}$, that has smallest volume of all simplexes
$T_b \subset \rr_+^n$ 
containing the both points $(\mu,0,1,\dots,1)$ and $(0,\nu,1,\dots,1)$ in their closure, 
and satisfying
$b_{1}=a_{1}$.
Then, certainly, $\vol T_{c} \le \vol T_{a}$.
To do that we need to minimize the function $V(b):=b_1\cdot\ldots\cdot b_n$
in the set 
\[
\{(b_1,\dots,b_n)\in\rr_+^{n}:\, b_1= a_{1}, 
\sum_{j=3}^\infty \frac 1{b_j}+\frac{\mu}{a_{1}}\le 1,
\sum_{j=3}^\infty \frac 1{b_j}+\frac{\nu}{b_{2}}\le 1\}.
\]
By standard calculations we obtain that the function $V$ attains the only minimum at the point 
$c = (c_{1},\dots,c_{n})$ where $c_{1}= a_{1}$,
$c_{2}= \frac{\nu}{\mu} a_{1}$, $c_{j}= (n-2)\frac{a_{1}}{a_{1}-\mu}$, $j=3,\dots,n$.
Therefore, we can estimate
\[
\begin{split}
\frac{\vol T_{a}}{\vol T} &\ge 
\frac{\vol T_{c}}{\vol T} 
=\frac{4x^2 \nu (n-2)^{n-2}a_{1}^n }{\mu n^n (a_{1}-\mu)^{n-2}}\\
&\overset{(\star)}{>}
\frac{4x^2 \nu (n-2)^{n-2}t^n }{\mu n^n (t-\mu)^{n-2}}
\underset{x\to 0}{\longrightarrow}
\frac{4(n-2)^{n-2}t^n}{n^n (t-1)^{n-2}} 
\overset{(\star)}{>} 1.
\end{split}
\]
Both inequalities marked with $(\star)$ hold because the function
$a\mapsto \frac{a^n}{(a-\theta)^{n-2}}$ is strictly increasing on 
$[\frac n2 \theta,+\infty)$---we use here inequalities $a_1>t>\frac n2$.
Thus, we get $\vol T_{a} > \vol T$ for sufficiently small $x>0$; a contradiction.
\end{proof}

\begin{remark}
The proof of Proposition~\ref{ex:three} is much simpler
if the family of pseudometrics $(\alpha_D)_D$ has the product property.
Indeed, by~{\rm(\ref{eq:prod})} we have
$\bb_{\alpha_{G_n}}(z) = \bb_{\alpha_{G_2}}(z_1,z_2) \times \ds^{n-2}$
for $z=(z_1,\dots,z_n)\in G_n$,
and from Proposition~\ref{thm:wprop}~(\ref{en:prod}) we get
\begin{multline*}
{\limsup_{0<x\to 0}\wu\alpha_{G_n}((x,0,\dots,0);(1,0,\dots,0))} \\
= \limsup_{0<x\to 0}\wu\alpha_{G_2}((x,0);(1,0)) 
\ge \sqrt 2 >1 = \wu\alpha_{G_n}(0;(1,0,\dots,0)),\\
\shoveleft{\limsup_{0<x\to 0}\twu\alpha_{G_n}((x,0,\dots,0);(1,0,\dots,0))}\\
\ge \frac {\sqrt 2}{\sqrt n} > \frac {1}{\sqrt{n-1}}
= \twu\alpha_{G_n}(0;(1,0,\dots,0)).
\end{multline*}
\end{remark}

\begin{proposition}
\label{ex:monotone}
There exist a domain $D\subset \cc^n$ $($for $n\ge 3)$ and an increasing sequence
of subdomains $D_m \nearrow D$ $(m\to\infty)$ such that $\twu\alpha_{D_m} \not\to \twu\alpha_D$
and $\wu\alpha_{D_m} \not\to \wu\alpha_D$.
\end{proposition}

\begin{proof}
Put $D:=G_n$ and fix a number $m\ge 1$. 
Let us consider two vectors $(1,0,1,\dots,1)$, $(0,m,1,\dots,1) \in \cc^n$.
As in the proof of Proposition~\ref{ex:three}, we show that the simplex
$T_m:=T_{(\frac n2, \frac {mn}2, n,\dots,n)}$ has smallest volume of all simplexes containing
the both vectors in their closure. Now, we take $D_m := G_n \cap \Psi^{-1}(T_m)$.
It is easy to see that $\bigcup_{m=1}^\infty D_m=D$ 
(for example, $G_n\cap\bb(0,\sqrt{\frac m2}) \subset D_m$).
Note that $D_m$ is a pseudoconvex Reinhardt domain. 
Therefore, we have 
$\bb_{\gamma_{D_m}}(0) = \conv D_m$ and $\bb_{\kappa_{D_m}}(0)=D_m$.
It implies that $\Psi\big(\bb_{\alpha_{D_m}}(0)\big) \subset T_m$, moreover,
$(1,0,1,\dots,1)$, $(0,\sqrt m,1,\dots,1) \in \partial\bb_{\alpha_{D_m}}(0)$.
Hence, we get from the minimality of $T_m$ that $\Psi\big(\bb_{\twu\alpha_{D_m}}(0) \big)=T_m$
(recall that $\bb_{\alpha_{D_m}}(0)$ is complete Reinhardt).
Consequently, we obtain $\twu\alpha_{D_m}\big(0;(1,0,\dots,0)\big) = \sqrt{\frac 2n}$,
$\wu\alpha_{D_m}\big(0;(1,0,\dots,0)\big) = \sqrt{2}$. In view of (\ref{eq:wugn}),
it finishes the proof.
\end{proof}

\begin{remark}
A similar counterexample 
is valid for dimension $2$ but only for $\wu\alpha$:  
if $D:= G_2$ and $D_m:= G_2 \cap \bb(0,m)$, then
$\wu\alpha_{D_m} \not\to \wu\alpha_D$.
We are not sure whether such a $2$--dimensional counterexample exists for $\twu\alpha$.
\end{remark}

\begin{proposition}\label{prop:cont}
Let $D$ be a domain in $\cc^2$. If $\eta_D\in\metrics(D)$ is a continuous pseudometric,
then $\twu \eta_D$ is upper semicontinuous.
\end{proposition}

\begin{corollary}\label{cor:gamma}
The pseudometric $\twu\gamma_D$ is upper semicontinuous for any domain $D\subset \cc^2$.
\end{corollary}

\begin{proof}[Proof of Proposition~\ref{prop:cont}]
Fix $z\in D$ and $X\in\cc^2$.

\emph{Case 1:
$\eta_D(z;Y)>0$ for any $Y\in\cc^2\setminus\{0\}$.}
Then $\eta_D$ is a metric in some neighborhood of the point $z$. 
The statement in Proposition~\ref{thm:wprop}~(\ref{en:cont}) is, in fact, local (see the proof)
and holds also for $\twu$.
Thus, $\twu\eta_D$ is continuous at the point $(z;X)$.

\emph{Case 2:
$\widehat\eta_D(z;Y)=0$ for some $Y\in\cc^2$.}
Therefore, the set of zeros of $\widehat\eta_D(z;\cdot)$ has the codimension at most $1$,
and consequently $\twu\eta_D(z;\cdot) = \widehat\eta_D(z;\cdot)$.
Recall that the pseudometric $\widehat\eta_D$
is upper semicontinuous (cf. \cite{bib:jp1}). Hence, we get
\[
\limsup_{(w,Y)\to (z,X)} \twu\eta_D(w;Y) 
\le \limsup_{(w,Y)\to (z,X)} \widehat\eta_D(w;Y) \le \widehat\eta_D(z;X) = \twu\eta_D(z;X).\qedhere
\]
\end{proof}

\begin{proposition}\label{prop:reinh}
Let $D \subset \cc^n$ be a pseudoconvex Reinhardt domain.
\begin{enumerate}
\item If $\gamma^{(k)}_D$ is a metric,
then $\twu\gamma^{(k)}_D$ and $\wu\gamma^{(k)}_D$ are continuous.
\item If $\kappa_D$ is a metric,
then $\twu\kappa_D$ and $\wu\kappa_D$ are continuous.
\item If $D$ is hyperconvex and $A_D$ is a metric,
then $\twu A_D$ and $\wu A_D$ are continuous.
\end{enumerate}
\end{proposition}

\begin{proof}
In view of Proposition~\ref{thm:wprop}~(\ref{en:cont}) it suffices to show
that the metrics $\gamma^{(k)}_D$, $\kappa_D$, and $A_D$ are continuous.

Hyperconvexity of $D$ immediately implies that $A_D$ is continuous (cf.\ \cite{bib:zwo2}).
If either $\gamma^{(k)}$ or $\kappa_D$ is a metric, then
the domain $D$ is pointwise $\kappa$--hyperbolic, and consequently Brody hyperbolic 
(i.e.\ all holomorphic maps from $\cc$ to $D$ are constant).
Due to the characterization theorem for hyperbolic pseudoconvex Reinhardt domains
(cf.\ \cite{bib:zwo1}, \cite{bib:zwo3}), $D$ is biholomorphic to a bounded domain 
(and so $\gamma$--hyperbolic) and taut.
Thus, the $\gamma$--hyperbolicity implies that $\gamma^{(k)}_D$ 
is a continuous metric (cf.\ \cite{bib:nik}).
The continuity of $\kappa_D$ follows from the tautness
(cf.\ \cite{bib:jp1}).
\end{proof}

\section{Formulae in elementary Reinhardt domains}

Let us introduce some notations concerning elementary Reinhardt domains.
We write $|z^\alpha|:= |z_1|^{\alpha_1}\dots |z_n|^{\alpha_n}$ 
for $\alpha=(\alpha_1,\dots,\alpha_n)\in\rr^n$ and $z\in\cc^n$, $z_j\neq 0$ if $\alpha_j<0$. 
For $\alpha \in (\rr\setminus\{0\})^n$  and $C>0$ define \emph{an elementary Reinhardt domain}
\[
D_{\alpha,C}:=\{ z\in\cc^n:\, |z^\alpha| <e^C 
\text{ and } \forall\,j=1,\dots,n:\ \alpha_j<0 \Rightarrow z_j\neq 0\}.
\]
We say that $D_{\alpha,C}$ is of \emph{rational type} if $\alpha\in\rr\cdot\zz^n$;
otherwise, it is of \emph{irrational type}.
Without loss of generality we may assume that $C=0$ and there exist $l\in\{0,\dots,n\}$
such that $\alpha_j<0$ for $j=1,\dots,l$ and $\alpha_j>0$ for $j=l+1,\dots,n$.
If $l<n$ then we put $t_l:= \min\{\alpha_{k+1}, \dots, \alpha_n\}$.
For $\alpha\in\zz^n$ and $r\in\nn$ put $\Phi(z):= z^\alpha$,
\[
\Phi_{(r)}(a)(X):= \sum_{\beta\in\zz_+^n,|\beta|=r} \frac {1}{\beta!}D^\beta\Phi(a)X^\beta,
\quad a\in D_\alpha, X\in\cc^n.
\]

The following formulae are known and collected in \cite{bib:jp2}.
\begin{proposition}
\label{prop:form}
Let $a\in D_\alpha$, $X\in\cc^n$.
Assume that $a_1\dots a_s\neq 0$, $a_{s+1}=\dots =a_n=0$ for some
$s\in\{l+1,\dots,n\}$. Put $r:= \alpha_{s+1}+\dots +\alpha_n$ if $s<n$ and $r:=1$ if $s=n$.
Consider the following four cases.
\begin{enumerate}
\item \label{en:acase}
$l<n$ and $D_\alpha$ is of rational type $($we may assume that $\alpha\in\zz^n$ and
$\alpha_1,\dots,\alpha_n$ are relatively prime$)$. Then:
\begin{align*}
\gamma_{D\alpha}(a;X)&= 
  \gamma_{\ds}\bigg(a^\alpha; a^\alpha\sum_{j=1}^n \frac{\alpha_jX_j}{a_j} \bigg),\\
A_{D\alpha}(a;X)&= \big(\gamma_\ds(a^\alpha; \Phi_{(r)}(a)(X)) \big)^{\frac 1r},\\
\kappa_{D\alpha}(a;X)&=
  \begin{cases}
    \gamma_\ds\Big((a^\alpha)^{\frac 1{t_l}}; (a^\alpha)^{\frac 1{t_l}} \frac 1{t_l}
      \sum_{j=1}^n \frac{\alpha_jX_j}{a_j} \Big) &\text{ if }s=n,\\
    \big(|a_1|^{\alpha_1}\dots |a_s|^{\alpha_s} 
         |X_{s+1}|^{\alpha_{s+1}}\dots |X_n|^{\alpha_n}\big)^{\frac 1r} &\text{ if } s<n.
  \end{cases}
\end{align*}
\item $l<n$ and $D_\alpha$ is of irrational type $($we may assume that $t_l=1)$. Then:
\begin{align*}
\gamma_{D\alpha}^{(k)}&\equiv 0, \quad k\ge 1, \\
A_{D\alpha}(a;X)&= 
  \begin{cases}
    0 &\text{ if }s=n,\\
    \big(|a_1|^{\alpha_1}\dots |a_s|^{\alpha_s} 
         |X_{s+1}|^{\alpha_{s+1}}\dots |X_n|^{\alpha_n}\big)^{\frac 1r} &\text{ if } s<n.
  \end{cases}\\
\kappa_{D\alpha}(a;X)&=
  \begin{cases}
    \gamma_\ds\Big(|a^\alpha|; |a^\alpha| \sum_{j=1}^n \frac{\alpha_jX_j}{a_j} \Big) 
      &\text{ if }s=n,\\
    \big(|a_1|^{\alpha_1}\dots |a_s|^{\alpha_s} 
         |X_{s+1}|^{\alpha_{s+1}}\dots |X_n|^{\alpha_n}\big)^{\frac 1r} &\text{ if } s<n.
  \end{cases}
\end{align*}
\item $l=n$ and $D_\alpha$ is of rational type $($we may assume that $\alpha\in\zz^n$ and
$\alpha_1,\dots,\alpha_n$ are relatively prime$)$. Then:
\begin{align*}
\gamma_{D_\alpha}^{(k)}(a;X) &= A_{D\alpha}(a;X) = 
  \gamma_{\ds}\bigg(a^\alpha; a^\alpha\sum_{j=1}^n \frac{\alpha_jX_j}{a_j} \bigg), \quad k\ge 1,\\
\kappa_{D\alpha}(a;X)&= 
  \kappa_{\ds\setminus\{0\}}\bigg(a^\alpha; a^\alpha\sum_{j=1}^n \frac{\alpha_jX_j}{a_j} \bigg).
\end{align*}
\item $l=n$ and $D_\alpha$ is of irrational type $($we may assume that $t_l=1)$. Then:
\begin{align*}
\gamma_{D_\alpha}^{(k)}(a;X) &= A_{D\alpha}(a;X) = 0, \quad k\ge 1,\\
\kappa_{D\alpha}(a;X)&= 
  \kappa_{\ds\setminus\{0\}}\bigg(|a^\alpha|; |a^\alpha|\sum_{j=1}^n \frac{\alpha_jX_j}{a_j}\bigg).
\end{align*}
\end{enumerate}
Moreover, if $\alpha\in\nn^n$ and $\alpha_1,\dots,\alpha_n$ are relatively prime, then
\[
\gamma_{D_\alpha}^{(k)}(a;X) = 
  \begin{cases}
    \big(\gamma_\ds(a^\alpha; \Phi_{(r)}(a)(X)) \big)^{\frac 1r} &\text{ if $r$ divides $k$},\\
    0, &\text{ otherwise },
  \end{cases}
\quad k\ge 1.
\]
\end{proposition}

\begin{proposition}
\label{prop:formulae}
Let the assumptions be the same as in Proposition~\ref{prop:form} and 
let $\eta\in\metrics(D)$ be any pseudometric such that $\eta \le \kappa_{D_\alpha}$.
Then:
\[
\wu\eta(a;X)=\twu\eta(a;X) =\widehat\eta(a;X).
\]
Moreover, if $\eta$ is one of $\gamma_{D_\alpha}^{(k)}, A_{D_\alpha}, \kappa_{D_\alpha}$,
then
\[
\widehat\eta(a;X)=
\begin{cases}
\eta(a;X)&\text{if }\  s\ge n-1,\\
0&\text{if }\  s<n-1.
\end{cases}
\]
\end{proposition}

\begin{proof}
One can see that the linear span of the set of zeros of
$\kappa_{D_\alpha}(a;\cdot)$
has the codimension either $0$ (if $s<n-1$) or $1$.
The same is true for $\eta$, and
therefore, the balls $\bb_{\wu\eta}(a)$, $\bb_{\twu\eta}(a)$,
and $\bb_{\widehat\eta}(a)$ obviously coincide.
When $\eta$ is one of the three above--mentioned metrics and $s\ge n-1$, the equality 
$\widehat\eta(a;X)=\eta(a;X)$ follows from the fomulae in Proposition~\ref{prop:form}.
\end{proof}

\end{document}